# Analysis and preconditioning of a probabilistic domain decomposition algorithm for elliptic boundary value problems

Francisco Bernal and Jorge Morón-Vidal

**Abstract**

PDDSparse is a new hybrid parallelisation scheme for solving large-scale elliptic boundary value problems on supercomputers, which can be described as a Feynman-Kac formula for domain decomposition. At its core lies a stochastic linear, sparse system for the solutions on the interfaces, whose entries are generated via Monte Carlo simulations. Assuming small statistical errors, we show that the random system matrix $\tilde{G}(\omega)$ is near a nonsingular M-matrix $G$, i.e. $\tilde{G}(\omega) + E = G$ where $\|E\|/\|G\|$ is small. Using nonstandard arguments, we bound $\|G^{-1}\|$ and the condition number of $G$, showing that both of them grow moderately with the degrees of freedom of the discretisation. Moreover, the truncated Neumann series of $G^{-1}$—which is straightforward to calculate—is the basis for an excellent preconditioner for $\tilde{G}(\omega)$. These findings are supported by numerical evidence.

**Keywords:** Domain decomposition, high performance computing (HPC), preconditioning, numerical PDE, Feynman-Kac formula, M-matrix, Neumann series, GMRES.

**MSC classifications:** 65F08, 65F50, 65N75, 65C05, 65Y05, 65D05

## 1 Introduction

The ongoing push of high performance parallel computing (HPPC) into the exascale has seen the emergence of truly massive supercomputers sporting millions of processors[1], frequently accelerated by unexpensive graphics processing units (GPUs). Such architectures are rather misaligned with state of the art domain decomposition (DD) algorithms [6, 14, 16, 17] for large-scale elliptic boundary value problems (BVPs), whose strong scalability is intrinsically capped above. The reason why is that subdomain-restricted BVPs are not well posed for lack of the solution along the fictitious interfaces, and thus require iterating towards them. This entails interprocessor communication and synchronisation, ultimately triggering Amdahl's law. As a rule of thumb, it is challenging to recruit more than thirty thousand processors in a DD simulation [2]. Another mismatch between hardware and software is that DD heavily relies on numerical linear algebra tasks, which are not always well suited to memory-limited GPUs.

---
[1] https://top500.org/lists/top500/2023/06/

Seeking to widen the strong scalability scope of DD, we proposed the PDDSparse algorithm in [3]. PDDSparse is an hybrid between Acebrón and Spigler's probabilistic domain decomposition framework [1], and linear algebra. The insight is to derive a Feynman-Kac formula for domain decomposition, i.e. a (sparse but unsymmetric) system $G\vec{u} = \vec{b}$ whose solution vector are the pointwise solutions of the original BVP on the fictitious interfaces (the skeleton, in DD nomenclature). Therefore, the system size is shrunk to $N = \mathcal{O}(\#DoF^{(D-1)/D})$, where $\#DoF$ is the number of degrees of freedom of the overall discretisation and $D$ the spatial dimension. The coefficients of the system are expectations of Feynman-Kac functionals of subdomain-bounded stochastic differential equations. They can be generated via Monte Carlo simulations—hence embarrassingly in parallel, and playing to the strength of GPUs.

The price to pay, though, is a distinctly more involved error structure with PDDSparse than with DD. Substituting sample estimators for the expected values in $G$ and $\vec{b}$ yields a (slightly) random system $\tilde{G}(\omega)\vec{u} = \vec{b}(\omega)$ instead, whose solution now carries statistical errors—on top to the usual discretisation error—resulting from the entanglement between noise and the spectral properties of the matrix through which noise is itself propagated. Preliminary experimental study had so far focused on scalability, error control, and mildly nonlinear elliptic equations [12], glossing over fundamental questions regarding the stability—and, in fact, the solvability—of the PDDSparse system. The aim of this paper is to rigorously explore the theoretical structure of PDDSparse.

We argue that $\tilde{G}(\omega)$ is a small random perturbation away from a diagonally dominant M-matrix $G$. By connectivity, the irreducibility of a certain principal minor of $G$ (and of $\tilde{G}(\omega)$) is first proved. Then, we motivate strong heuristics to claim that the off-diagonal terms of $G$ are asymptotically nonpositive. Interpolatory and Feynman-Kac arguments lead to $G$ being row-diagonally dominant, whence invertibility and (later) M-matrixness of $G$ follow in a standard way.

Two consequences are the nonnegativity of $G^{-1}$ and the convergence of the Neumann ( matrix geometric) series of $G^{-1}$. We use the former to bound $\|G^{-1}\|_\infty$ in terms of expected first-exit and boundary residence times of associated SDEs. The proof is, in our opinion, interesting in itself, since it combines stochastic calculus and matrix analysis. Moreover, the bound is geometrically interpretable and allows us to bound the condition number of $G$ with respect to $N$ in several typical scenarios. Asymptotically, $\kappa_\infty(G)$ grows as $\mathcal{O}(N)$ at most (i.e. roughly with $(\#DoF)^{1/2}$), thus confirming the stability of the PDDSparse system.

Regarding the Neumann polynomial for $G^{-1}$ (truncated at $t$), we employ it for a simple, tailored left preconditioner $P_t^{-1}$ for the actual (noisy) matrix $\tilde{G}(\omega)$. While it reduces both the condition number of $P_t^{-1}\tilde{G}$ and the number of GMRES iterations, the wallclock time savings are modest, if any. In order to enhance the spectrum of $P_t^{-1}\tilde{G}$ while keeping $t$ low, we have borrowed an idea from the recent paper [18] by Zheng, Xi, and Saad. Their insight is to include spectral information related to the correction $P_t^{-1} - \tilde{G}^{-1}$ via an incomplete Arnoldi procedure—which is particularly suitable for large, sparse matrices as $\tilde{G}$. It turns out that $P_t^{-1}$ works best within that framework—although the latter is designed for quite another situation. In our largest example (with $\tilde{G}$ of almost order $10^5$), the spectrum of the preconditioned matrix was clusterised by orders of magnitude, and the GMRES iterations collapsed from 193 to just 2.

The remainder of the paper is organised as follows. Section 2 succintly reviews the PDDSparse algorithm. (However, the extension to reflecting BCs is a novelty.) Section 3 motivates and states the main approximation. Based on it, the properties of a "limit" matrix $G$ (noise-free), of $\|G^{-1}\|_\infty$, and of $\kappa_\infty(G)$ are rigorously derived in Section 4. The structure of $G$ is exploited in Section 5 to design an approximation to $\tilde{G}^{-1}(\omega)$ with good preconditioning properties. The theory is empirically validated in Section 6, and Section 7 concludes the paper.

## 2  The PDDSparse algorithm

This section is abridged from [3], except for the extension to reflecting BCs.

### 2.1  Feynman-Kac formula for elliptic BVPs with mixed BCs

Let $\Omega$ be a bounded, open, connected set in $\mathbb{R}^D$ such that $\partial\Omega = \partial\Omega_A \cup \partial\Omega_R$ and $\partial\Omega_A \cap \partial\Omega_R = \emptyset \neq \partial\Omega_A$. ($D \in \mathbb{N}$ in general, but in this paper, $D = 2$.) Let

$$\begin{cases} \overbrace{\frac{1}{2}\sum_{1\leq i\leq j\leq D} a_{ij}(\mathbf{x})\frac{\partial^2 u(\mathbf{x})}{\partial x_i \partial x_j} + \sum_{1\leq k\leq D} b_k(\mathbf{x})\frac{\partial u(\mathbf{x})}{\partial x_k}}^{:=\mathcal{L}u} + c(\mathbf{x})u(\mathbf{x}) + f(\mathbf{x}) = 0 & \text{in } \Omega, \\ u(\mathbf{x}) = g(\mathbf{x}) & \text{on } \partial\Omega_A, \\ \mathbf{N}^\top(\mathbf{x})\nabla u(\mathbf{x}) = \varphi(\mathbf{x})u(\mathbf{x}) + \psi(\mathbf{x}) & \text{on } \partial\Omega_R, \end{cases} \quad (1)$$

be an elliptic BVP with Dirichlet (for "absorbing") BCs on $\partial\Omega_A$ and Neumann/ Robin (for "reflecting") boundary conditions (BCs) on $\partial\Omega_R$, where $\mathbf{N}$ is the outward normal to $\partial\Omega$. We shall assume that $c, \varphi \leq 0$ and that a unique, smooth solution to (1) exists. (See [8] for the required conditions on $\partial\Omega$ and the coefficients.) Consider the system of stochastic differential equations (SDEs):

$$\begin{cases} d\mathbf{X}_t = \mathbf{b}(\mathbf{X}_t)\,dt + \sigma(\mathbf{X}_t)\,d\mathbf{W}_t - \mathbf{N}(\mathbf{X}_t)\,d\xi_t & \mathbf{X}_0 = \mathbf{x}_0 \in \Omega, \\ dY_t = c(\mathbf{X}_t)Y_t\,dt + \varphi(\mathbf{X}_t)Y_t\,d\xi_t & Y_0 = 1, \\ dZ_t = f(\mathbf{X}_t)Y_t\,dt + \psi(\mathbf{X}_t)Y_t\,d\xi_t & Z_0 = 0, \\ d\xi_t = \mathbb{1}_{\{\mathbf{X}_t \in \partial\Omega_R\}}\,dt & \xi_0 = 0, \end{cases} \quad (2)$$

where $\mathbf{W}_t$ is the $D$–dimensional Wiener process and $[a_{ij}] = \sigma\sigma^\top$ is positive definite. It is helpful to think of each realisation of the bounded process $(\mathbf{X}_t)_{0\leq t\leq \tau^\Omega} \in$ $\overline{\Omega}$ as a stochastic "trajectory" which starts off at $\mathbf{x}_0$, diffuses across $\Omega$, is normally reflected on $\partial\Omega_R$ (assuming $\partial\Omega_R \neq \emptyset$) if it hits it, and is finally absorbed upon hitting $\partial\Omega_A$. The random variables $\tau^\Omega(\mathbf{x}_0) = \min_{t>0}\{\mathbf{X}_t \in \partial\Omega_A \,|\, \mathbf{X}_0 = \mathbf{x}_0\} > 0$ and $\xi_\tau \geq 0$ are the first exit time (FET) and local time of the trajectory, respectively. (Intuitively, $\xi_\tau$ is the time spent on $\partial\Omega_R$.)

Then, under somewhat more restrictive assumptions over the $\{a_{ij}\}, \{b_k\}$, and $\partial\Omega$ (see [5] for details), and provided that $\mathbb{E}[\tau^\Omega(\mathbf{x}_0)] < \infty$ if $c \equiv 0$, it holds [10]

$$u(\mathbf{x}_0) = \mathbb{E}[g(\mathbf{X}_{\tau^\Omega})Y_{\tau^\Omega} + Z_{\tau^\Omega} \,|\, \mathbf{X}_0 = \mathbf{x}_0]. \qquad (3)$$

The following facts will later be used in the analysis:

- The "trajectory score" $g(\mathbf{X}_{\tau^\Omega})Y_{\tau^\Omega} + Z_{\tau^\Omega}$ depends on the coefficients $c$ and $f$ along the trajectory; on the value of the Dirichlet BC at the first exit point $\mathbf{X}_\tau \in \partial\Omega_A$; and on the coefficients $\varphi, \psi$ on the reflections on $\partial\Omega_R$ (if any); but not explicitly on the first and second derivative coefficients in (1).

- On the other hand, given a starting point $\mathbf{x}_0$, the trajectories across $\Omega$ are determined solely by $\mathcal{L}$ and the normal to $\partial\Omega_R$ (if any), but not $c, f, g, \varphi, \psi$.

- $0 < Y_{\tau^\Omega} \leq 1$ depends only on $c$ and $\varphi$; and $\mathbb{E}[Y_{\tau^\Omega}] = 1$ iff $c \equiv \varphi \equiv 0$. Also,

$$\mathbb{E}[g(\mathbf{X}_{\tau^\Omega})Y_{\tau^\Omega}] \leq \mathbb{E}[g(\mathbf{X}_{\tau^\Omega})]. \qquad (4)$$

- Setting $c \equiv 0$, $f \equiv 1$, the expected FET is the solution of the Poisson BVP

$$\mathcal{L}\mathbb{E}[\tau^\Omega(\mathbf{x})] + 1 = 0 \text{ in } \Omega, \; \mathbb{E}[\tau^\Omega(\mathbf{x})] = 0 \text{ on } \partial\Omega_A, \; \mathbf{N}^\top \nabla\mathbb{E}[\tau^\Omega(\mathbf{x})] = 0 \text{ on } \partial\Omega_R. \quad (5)$$

## 2.2 Discretisation

We shall introduce the discretisation elements of PDDSparse in an intuitive way, mostly referring to Figure 1. The reader is directed to [3] for details.

For the purpose of domain decomposition, the bidimensional domain $\Omega$ is embedded into a grid of squares, whereby two types of nonoverlapping subdomains are defined. Squares fully immersed in $\Omega$ (save perhaps isolated points on $\partial\Omega$) give rise to *floating* subdomains. Otherwise, they give rise to *perimeter* subdomains. We call the set of squares intersecting $\Omega$ the *tessellation*. The underlying grid lines are segmented into *interfaces* (the subdomain fictitious boundaries), and discretised into $N$ *knots*[2] on which the *nodal solutions* $u_1 \approx u(\mathbf{x}_1), \ldots, u_N \approx u(\mathbf{x}_N)$ will be calculated. Knots are equispaced according to a constant *linear density of knots* $n = 1/\Delta z$, in such a way that there is always a knot at the grid crossings and at the intersection with the actual boundary $\partial\Omega$. Every knot $\mathbf{x}_i$ has a nodal *patch* $\Pi(\mathbf{x}_i)$ associated to it, defined as the union of all the subdomains that $\mathbf{x}_i$ belongs to. *Floating patches* are composed of floating subdomains only, and are either rectangular (by default) or square (if they are associated to a crossing knot). Otherwise, they are a *perimeter patch*. Unlike floating patches, the latter may be irregularly shaped, and are exposed to the actual BCs on at least one of the four "sides" of the patch. (They are denoted as $\{E, N, W, S\}$ standing for East, North, etc.) With each interface $O \in \{E, N, W, S\}$ of patch[3] $\Pi(\mathbf{x}_i)$, namely $\partial^O \Pi(\mathbf{x}_i)$, we associate a *patch side stencil* $\mathcal{S}_i^O$ made up of the knots sitting on the patch interface (possibly *elongated* beyond either or both patch corners on that interface, with no gaps). The idea is that $u|_{\partial \Pi^O(\mathbf{x}_i)}$ (the solution on the given interface) is to be interpolated over the knots in its stencil and their (yet unknown) nodal values.

---

[2] We reserve the word "nodes" for the parallel computer nodes.

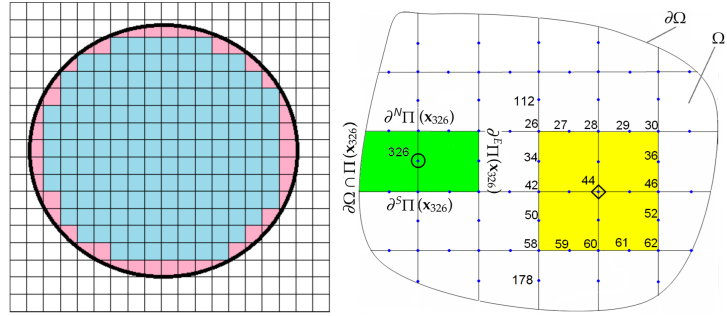

Figure 1: Left: a circular domain $\Omega$ discretised into a PDDSparse tessellation of perimeter (pink) and floating (blue) subdomains. Right: a toy discretisation. The three interfaces of the perimeter patch $\Pi(\mathbf{x}_{326})$ (green) are labeled. Also, the knots on $\partial \Pi(\mathbf{x}_{44})$, and the (elongated) stencil $\mathcal{S}_{44}^W = \{112, 26, 34, 42, 50, 58, 178\}$.

### 2.3 The PDDSparse linear system

The insight is to derive a "Feynman-Kac formula for domain decomposition". By ellipticity of $\mathcal{L}$, the solution $u|_{\Pi(\mathbf{x}_i)}$ inside any patch $\Pi(\mathbf{x}_i)$ for $1 \leq i \leq N$ is also the solution of the patch-restricted BVP:

$$\begin{cases} \mathcal{L}u + cu + f = 0 & \text{if } \mathbf{x} \in \Pi(\mathbf{x}_i), \\ u = u|_{\partial \Pi(\mathbf{x}_i)} & \text{if } \partial^O \Pi(\mathbf{x}_i) \neq \emptyset \text{ and } \mathbf{x} \in \partial^O \Pi(\mathbf{x}_i) \text{ for } O = \{E, N, W, S\}, \\ u = g & \text{if } \Pi(\mathbf{x}_i) \cap \partial \Omega_A \neq \emptyset \text{ and } \mathbf{x} \in \Pi(\mathbf{x}_i) \cap \partial \Omega_A, \\ \mathbf{N}^\top \nabla u = \varphi u + \psi & \text{if } \Pi(\mathbf{x}_i) \cap \partial \Omega_R \neq \emptyset \text{ and } \mathbf{x} \in \Pi(\mathbf{x}_i) \cap \partial \Omega_R. \end{cases} \quad (6)$$

(Note that the last two equations only apply on perimeter patches.) The nodal solution $u(\mathbf{x}_i)$ can then be written in terms of the Feynman-Kac formula as

$$u(\mathbf{x}_i) = \mathbb{E}\left[Z_{\tau_i^\#} \,\middle|\, \mathbf{X}_0 = \mathbf{x}_i\right] + \mathbb{E}\left[Y_{\tau_i^\#} g\left(\mathbf{X}_{\tau_i^\#}\right) \mathbb{1}_{\{\mathbf{X}_{\tau_i^\#} \in \partial \Omega_A\}} \,\middle|\, \mathbf{X}_0 = \mathbf{x}_i\right]$$
$$+ \sum_{O \in \{E,N,W,S\}} \mathbb{1}_{\partial^O \Pi(\mathbf{x}_i) \neq \emptyset} \mathbb{E}\left[Y_{\tau_i^\#} \mathbb{1}_{\{\mathbf{X}_{\tau_i^\#} \in \partial^O \Pi(\mathbf{x}_i)\}} u|_{\partial^O \Pi(\mathbf{x}_i)}(\mathbf{X}_{\tau_i^\#}) \,\middle|\, \mathbf{X}_0 = \mathbf{x}_i\right], \quad (7)$$

---

[3]In general, each patch will be associated to multiple knots.

where the diffusions $(\mathbf{X}_t)_{0 \le t \le \tau_i^\#}$ are now confined to the patch $\Pi(\mathbf{x}_i)$. They start at knot $\mathbf{x}_i$, may be reflected on the actual reflecting boundary iff the patch is a perimeter one facing $\partial \Omega_R$, and finish either on the patch interfaces, or on $\partial \Omega_A$ (where they read the actual Dirichlet BC $g$) in case of perimeter patches exposed to it. We define $\tau_i^\# = \tau^{\Pi(\mathbf{x}_i)}(\mathbf{x}_i)$ as the *first exit time* (FET) from patch $\Pi(\mathbf{x}_i)$ starting from $\mathbf{x}_i$. The local time $\xi_i^\#$ (affecting $Y_{\tau_i^\#}$ and $Z_{\tau_i^\#}$) is defined likewise.

The above expression is purely formal, since the interfacial BCs $u|_{\partial^O \Pi(\mathbf{x}_i)}$ are unknown. Approximating them by linear interpolation, the Dirichlet BCs are

$$u|_{\partial^O \Pi(\mathbf{x}_i)}(\mathbf{x}) \approx \sum_{\mathbf{x}_j \in \mathcal{S}_i^O} u(\mathbf{x}_j) H_{ij}^O(z) \qquad \text{if } \mathbf{x} \in \partial^O \Pi(\mathbf{x}_i). \tag{8}$$

$H_{ij}^O(z)$ is the $j^{th}$ interpolatory cardinal function in the segment defined by the collinear set of knots in $\mathcal{S}_i^O$, and parameterised by $z \in [z_{min}, z_{max}]$. Then,

$$u(\mathbf{x}_i) \approx \mathbb{E}\left[Z_{\tau_i^\#} \,\Big|\, \mathbf{X}_0 = \mathbf{x}_i\right] + \mathbb{E}\left[Y_{\tau_i^\#} g\left(\mathbf{X}_{\tau_i^\#}\right) \mathbb{1}_{\{\mathbf{X}_{\tau_i^\#} \in \partial \Omega_A\}} \,\Big|\, \mathbf{X}_0 = \mathbf{x}_i\right]$$
$$+ \sum_{O \in \{E,N,W,S\}} \mathbb{1}_{\{\partial^O \Pi(\mathbf{x}_i) \ne \emptyset\}} \sum_{\mathbf{x}_j \in \mathcal{S}_i^O} \mathbb{E}\left[Y_{\tau_i^\#} \mathbb{1}_{\{\mathbf{X}_{\tau_i^\#} \in \partial^O \Pi(\mathbf{x}_i)\}} H_{ij}^O\!\left(z(\mathbf{X}_{\tau_i^\#})\right) \,\Big|\, \mathbf{X}_0 = \mathbf{x}_i\right] u(\mathbf{x}_j). \tag{9}$$

Denoting as $u_k \approx u(\mathbf{x}_k)$ ($k = 1, \ldots, N$) the resulting approximate nodal solution (with interpolation error), and rearranging terms, the system arises:

$$G(u_1, \ldots, u_N)^\top = \vec{b}, \qquad \text{or} \qquad G\vec{u} = \vec{b}, \tag{10}$$

whose entries are (writting $H_{ij}^O(\mathbf{X}_{\tau_i^\#})$ instead of $H_{ij}^O[z(\mathbf{X}_{\tau_i^\#})]$ for brevity)

$$G_{ij} = \begin{cases} 1 & \text{if } i = j, \\ -\mathbb{E}\left[Y_{\tau_i^\#} \mathbb{1}_{\{\mathbf{X}_{\tau_i^\#} \in \partial^O \Pi(\mathbf{x}_i)\}} H_{ij}^O(\mathbf{X}_{\tau_i^\#}) \,\Big|\, \mathbf{X}_0 = \mathbf{x}_i\right] & \text{if } \mathbf{x}_j \in \mathcal{S}_i^O,\, O \in \{E, N, W, S\}, \\ 0 & \text{otherwise}, \end{cases} \tag{11}$$

and

$$(\vec{b})_i = \mathbb{E}\left[Z_{\tau_i^\#} \,\Big|\, \mathbf{X}_0 = \mathbf{x}_i\right] + \mathbb{E}\left[Y_{\tau_i^\#} g\left(\mathbf{X}_{\tau_i^\#}\right) \mathbb{1}_{\{\mathbf{X}_{\tau_i^\#} \in \partial \Omega_A\}} \,\Big|\, \mathbf{X}_0 = \mathbf{x}_i\right]. \tag{12}$$

For radial basis function (RBF) interpolation, if the stencil of $\partial^O \Pi(\mathbf{x}_i)$ is given by knots at positions $z_1, \ldots, z_m$ along a line, then

$$H_{ij}^O(z) = \sum_{k=1}^m \left(\Phi^{(i)}\right)^{-1}_{km} \phi_c(|z - z_k|), \tag{13}$$

where $\phi_c$ is the RBF with shape parameter $c > 0$—in this paper we have used the Gaussian, $\phi_c(z) = e^{-z^2/c^2}$—and $[\Phi^{(i)}]_{ln} = \phi_c(|z_l - z_n|)$ is the interpolation matrix.

By interpolating $f(z) = 1$ in the stencil interval, it is straightforward that

$$\sum_{1 \le j \le m} H_{ij}^O(z) \approx 1 \tag{14}$$

over the interval. For the Gaussian RBF, the difference goes to zero exponentially fast as $\Delta z \to 0$, $c \to \infty$ ( $\Delta z$ being the constant separation between RBFs).

**Remark 1** *In* (13), *note that the stencils need not be similar, or have the same spacing, on all interfaces. If fact, patches need not be rectangular and might even overlap. This allows, in principle, for local refinement, or curved interfaces close to the boundaries.*

## 3 Main approximation

Understanding the interplay between stochastics, interpolation, and matrix algebra, is the key to analysing PDDSparse. Figure 2 will help build intuition. A perimeter patch corresponding to a knot sitting on an interfacial crossing (say $\mathbf{x}_i$) is depicted[4]. Two stochastic trajectories stemming from it are sketched. One exits the patch at point $\mathbf{X}^{(1)}_{\tau_i^\#}$ of one interface (say $\partial \Pi^N(\mathbf{x}_i)$), where it reads the values at $\mathbf{X}^{(1)}_{\tau_i^\#}$ of *all of* the cardinal functions $H^N_{ij}$ associated with the interpolation stencil of that interface. The nonzero off-diagonal matrix entries on row $i^{th}$ of $G$ are the $Y_{\tau_i^\#}$-weighted average over all such contributions. Because this particular patch is a perimeter one, a portion of trajectories will hit the actual domain boundary $\partial \Omega$, where they read the Dirichlet BC, and contribute the right-hand side vector entry $b_i$ instead. If the BCs on that portion of $\partial \Omega$ were of Neumann or Robin type[5], the trajectory would be reflected back into $\Pi(\mathbf{x}_i)$, and eventually hit an interface.

In Figure 2, one of the $H^N_{ij}(z)$ functions defined over $\partial^N \Pi(\mathbf{x}_i)$ is depicted. Not by chance, it resembles a scaled sinc function. In fact, for many RBFs the cardinal function of equispaced RBF interpolation over $[-a, a]$ tends to $\text{sinc}(z/\Delta z)$ as $a \to \infty$ and the shape parameter $c \to \infty$ [4]. But when the interval covered by the stencil is finite, the RBF cardinal functions cannot be approximated by a sinc near the interval ends (see Figure 2, right). Therefore, it makes sense to *elongate* the patch side stencils beyond the patch corners (including knots from the previous and next interfaces along the same line), so that $H^N_{ij} \approx \text{sinc}((z - z_j)/\Delta z)$ *along the patch interface* $\partial^N \Pi(\mathbf{x}_i)$—*which is where trajectories do impact.* Indeed, the other ingredient of $G_{ij}$ is the *probability densities of exit points* along the interfaces of the patch boundary, namely $\rho_i$ and $\{\rho_i^O\}$:

$$\oint_{\partial \Pi(\mathbf{x}_i)} \rho_i(\mathbf{z})\, d\ell := \sum_{O \in \{E,N,W,S\}} \mathbb{P}[\mathbf{X}_{\tau_i^\#} \in \partial^O \Pi(\mathbf{x}_i) \mid \mathbf{X}_0 = \mathbf{x}_i] \ =: \sum_{O \in \{E,N,W,S\}} \int_0^{|\partial^O \Pi(\mathbf{x}_i)|} \rho_i^O(z)\, dz,$$

where $|\partial^O \Pi(\mathbf{x}_i)|$ is the length on that patch interface. Since some trajectories will be absorbed on $\partial \Omega_A$ if the patch is exposed to it, $\partial \Omega_A$ works as a "leak" of

---

[4]The arguments are the same for all patches regardless of their shape.
[5]To simplify the discussion, we assume that no patch is facing both Dirichlet and reflecting BCs.

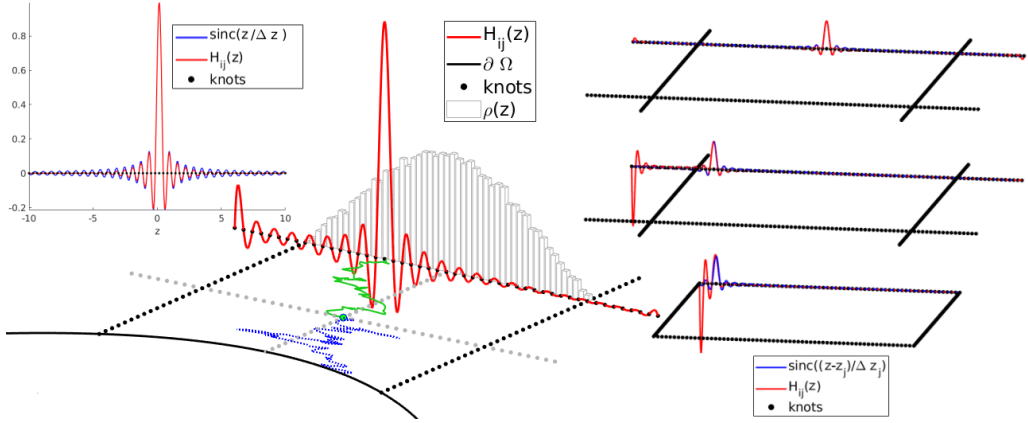

Figure 2: Left: a floating patch associated to a crossing knot, from which two exemplary patch-confined diffusions start. The blue (dotted) one hits the boundary, while the green one is shown hitting on exiting the patch one of the cardinal functions on that patch interface. Underlaid is an histogram of first-exit points on that side (not at scale). As shown in the top left inset, the cardinal function and the sinc are very similar far from the stencil ends. Right, bottom to top: without elongated stencils (bottom), the Gaussian RBF cardinal function departs visibly from the sinc close to the corners, which then coincide with the stencil ends. This is fixed by elongating the stencils beyond the corners.

first exit points on the interfaces on perimeter subdomains, and

$$\oint_{\partial \Pi(\mathbf{x}_i)} \rho_i(\mathbf{z}) \, d\ell \begin{cases} = 1, & \text{for floating, or perimeter patch facing } \partial \Omega_R \\ < 1, & \text{for perimeter patch facing } \partial \Omega_A. \end{cases} \quad (15)$$

From (11), the link between algebra, interpolation and stochastics is, then:

$$G_{ij} = - \sum_{O \in \{E,N,W,S\}} \int_0^{|\partial^O \Pi(\mathbf{x}_i)|} H_{ij}^O(z) \rho_i^O(z) \, dz. \quad (16)$$

In Figure 2, an actual patch side density $\rho_i^O$ is underlaid over an actual $H_{ij}^N(z)$ cardinal function. The main "lobe" makes up about 60% of the total area. Since the area of the positive lobes is significantly larger than that of the negative lobes[6], and the oscillation period of the sinc (i.e. the internodal separation $\Delta z \ll |\partial^O \Pi(\mathbf{x}_i)|$) is much shorter than the typical length of the density $\rho_i^O$, the interfacial integrals in (16) are intuitively expected to be positive. In order to occur otherwise, $\rho_i^O(z)$ would need to have very narrow peaks coupled with the negative lobes of the cardinal functions. This can not happen with smooth enough[7] first and second derivative coefficients in the differential generator $\mathcal{L}$. Moreover, the aforementioned "coupling" could be fixed by refining the discretisation—by just decreasing $\Delta z$ so that the peaks fall on positive lobes.

---

[6]The difference is about one, since $\int_{-\infty}^{+\infty} \text{sinc}(z) \, dz = 1$.

In sum, we have motivated:

**Heuristic (H).** With elongated stencils, $c$ large, and $\Delta z$ small enough, the cardinal functions $H_{ij}^O(z)$ of Gaussian-RBF interpolation are approximately $\mathrm{sinc}\left(\frac{z-z_j}{\Delta z}\right)$.

**Remark 2** *Interfaces terminating on $\partial\Omega$ cannot be elongated—whence finite end-of-interval oscillations such as in Figure 2 (bottom right) will appear.[8] An idea would be to elongate those stencils with "ghost knots" outside $\Omega$; we leave this for future work.[9]*

## 4 Analysis of $G$

### 4.1 Irreducibility of $G^*$

The sparsity pattern of $G$ reflects the knot stencils. Knot $i$ is connected to knot $j$ by a directed edge $i \to j$ iff $G_{ij} \neq 0$. By (11), the filled columns $\{j\}$ of row $i$ of $G$ are the knots $\{\mathbf{x}_j\}$ in the stencil of $\mathbf{x}_i$, plus the diagonal. The exception are the $m$ knots on $\partial\Omega_A$ (whose value is the Dirichlet BC), connected to themselves only. Ordering knots so that they are indexed last, $G$ is structured as[10]

$$G = \begin{bmatrix} G^* & [*] \\ [0] & I_m \end{bmatrix}, \qquad (17)$$

where $I_m$ is the identity of order $m$, $[*]$ has size $(N-m) \times m$, $[0]$ is a matrix of zeros, and $G^*$ is the $(N-m) \times (N-m)$ principal minor. Since $\det(G) = \det(G^*)$,

**Lemma 1** *$G$ is invertible iff $G^*$ is invertible.*

Thus, $G^*$ is the matrix (and the associated digraph) obtained by deleting the Dirichlet BC knots and their interconnections with other knots. Due to the identity block, the digraph (directed graph) of $G$ cannot be strongly connected—unlike, as we shall prove, the digraph of $G^*$. A preparatory lemma comes next, whereby knots sitting on patch interfaces are classified into "end knots" (those sitting at crossings or on the boundary $\partial\Omega$) and "mid knots". For a square, the following directed edges exist connecting its interfacial knots.

**Lemma 2** *Any PDDSparse patch boundary has at least one interface—let it be W. i) Mid knots of W are connected to mid knots of N and S, and vice versa (a.v.v.). ii) Mid knots of W are connected to all knots of E, a.v.v. iii) SW end knot is connected to all knots of N and E, a.v.v. iv) NW end knot is connected to all knots of S and E, a.v.v.*

---

[7]Not regular enough coefficients render the Feynman-Kac formula (3) void, in the first place.
[8]Indeed, the $G_{ij}$ of next-to-boundary knots are experimentally the most likely to overshoot zero.
[9]Preliminary tests do dampen oscillations and overshoots, but make little difference in practice.
[10]Neither here nor anywhere else do asterisks stand for complex conjugate.

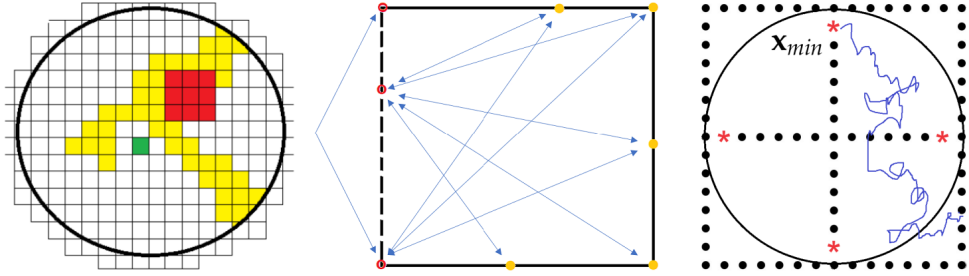

Figure 3: Left: constructing the tessellation in Figure 1 according to Theorem 4—a possible initial (red) and valid intermediate (red and yellow) tessellations. Centre: the directed edges between one side of a square ($W$) and the other three. Right: the geometrical argument employed for deriving (34).

*Proof.* The directed edges quoted above exist iff one can draw a Feynman-Kac trajectory inside the patch of the starting knot linking the starting knot and the interface of the absorption knot. They are depicted in Figure 3 (centre). □

**Definition 3** *Let $\mathcal{T}_M$ be a PDDSparse tessellation with M subdomains. We call it valid if: i) $M \geq 9$; ii) it contains a subtessellation $\mathcal{T}_9 \subset \mathcal{T}_M$ of 3×3 floating subdomains; and iii) every subdomain lies side by side with at least another one.*

For instance, the yellow tessellation $\mathcal{T}_{50}$ in Figure 3 (left) is not valid, but the union of the yellow and the red ($\mathcal{T}_9$) is ($\mathcal{T}_{59} := \mathcal{T}_{50} \cup \mathcal{T}_9$). On the other hand, $\mathcal{T}_{59}$ plus the green square is not a valid tessellation, because the green subdomain is not connected sideways with any subdomain of $\mathcal{T}_{59}$ (only at a corner).

**Theorem 4** *Let $\mathcal{T}_M$ be a valid PDDSparse tessellation. Then $G^*$ is irreducible.*

*Proof.* First, we delete from $G^*$ the edges induced by elongated stencils: let

$$\hat{G}^*_{ij} = \begin{cases} G^*_{ij} & \text{if } \mathbf{x}_j \in \partial \Pi(\mathbf{x}_i) \text{ (belongs to the stencil of } \mathbf{x}_i \text{ on its patch),} \\ 0 & \text{if not (removing elongation edges from the digraph).} \end{cases} \quad (18)$$

With every PDDSparse tessellation of $n$ subdomains, denoted as $\mathcal{T}_n$, let us associate a digraph $\mathcal{D}_n$ with adjacency matrix $A_n$. Figure 1 (left) shows a PDDSparse tessellation with $n = 203$ subdomains, $\mathcal{T}_{203}$. In Figure 3 (left), a 3×3 subtessellation $\mathcal{T}_9$—as assumed by the theorem—is highlighted in red. Clearly, it is possible to reach any final valid PDDSparse tessellation $\mathcal{T}_M$ by attaching $M - 9$ subdomains by the side, one at a time, starting from $\mathcal{T}_9$. A possible intermediate tessellation $\mathcal{T}_{50}$ to $\mathcal{T}_{203}$ is highlighted in yellow in Figure 3 (left).

We shall prove by induction that, if the digraph $\mathcal{D}_n$ associated to the knots of a tessellation $\mathcal{T}_n$ is strongly connected, then the digraph of the tessellation $\mathcal{T}_{n+1}$ defined by attaching a new square by the side to $\mathcal{T}_n$ is also strongly connected.

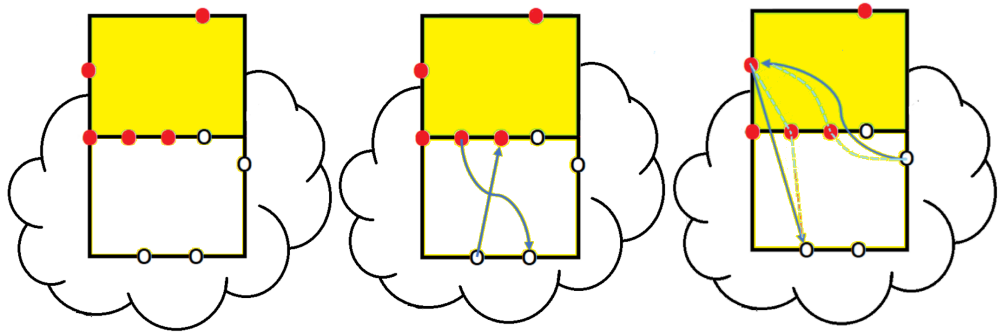

Figure 4: Attaching a new square to an extant tessellation (clouds) according to Theorem 4. The new (green) and at least one "old" (yellow) subdomain share one interface (left), whose knots are connected to the old digraph (centre). Then, the rest of new knots are also connected through the common interface (right).

To verify this, the reader is referred to Figure 4. Attaching a new subdomain $\Omega_{n+1}$ means that there emerges at least one interface, $I(n+1, r)$, between a subdomain $\Omega_r$ ($r \leq n$) belonging to the old tessellation $\mathcal{T}_n$ and the new subdomain. Two kinds of new knots arise: those sitting on the interface $I(n+1, r)$ (necessarily including new mid knots), and at least on another side—different from $I(n+1, r)$—of the new subdomain $\Omega_{n+1}$. Consider first one new mid knot on $I(n+1, r)$. By Lemma 2, it is connected to at least one old knot of $\Omega_r$; and at least one old knot of $\Omega_r$ is connected to the new mid knot. Therefore, all mid knots on the shared new interface $I(n+1, r)$ are strongly connected to $\mathcal{T}_n$.

Again by Lemma 2, all of the new knots (which are not mid knots of $I(n+1,r)$) on the new subdomain are connected to the mid knots of $I(n+1, r)$, and vice versa. And, through the latter, the former are connected to the old knots in $\mathcal{T}_n$, and vice versa. Consequently, the old digraph $\mathcal{D}_n$ is augmented by all the new knots and edges in such a way that the resulting new digraph $\mathcal{D}_{n+1}$ remains strongly connected. Equivalently, the adjacency matrix $A_n$ is enlarged (with as many rows and columns as there are new knots) and filled into $A_{n+1}$ in such a way that $A_{n+1}$ is irreducible too. The base case of the induction is $\mathcal{T}_9$, which can easily be shown to be irreducible.

Summing up, there is a constructive procedure which produces the tessellation $\mathcal{T}_M$ and an associated adjacency matrix $A_M$ identical to that of the $\hat{G}^*$ matrix of any valid PDDSparse discretisation *without elongated stencils*, and which is irreducible. To avoid a contradiction, $\hat{G}^*$ must be irreducible, too.

Finally, if $\hat{G}^*$ is irreducible, then $G^*$ is irreducible too, since adding edges to a strongly connected graph preserves the strong connectivity property. □

## 4.2 $G$ is an invertible M-matrix under (H)

**Theorem 5** *Under the proviso of heuristic* **(H)** $G_{ij} \leq 0$ *if* $i \neq j$ *as* $\Delta z \to 0^+$.

*Proof:* By **(H)**, $H_{ij}^O(z) \to \text{sinc}\left(\frac{z-z_j}{\Delta z}\right)$ as $\Delta z \to 0^+$. The normalised $\text{sinc}(x)$ is

$$\text{sinc}(x) = \sin(\pi x)/(\pi x) \text{ if } x \neq 0, \qquad \text{sinc}(0) = 1. \tag{19}$$

(The name reflects that $\int_{-\infty}^{+\infty} \text{sinc}(x)\,dx = 1$.) It has the "Kronecker delta" property: if $m, n \in \mathbb{Z}$, then $\text{sinc}(n) = \delta_{n0}$ (where $\delta_{nm} = 1$ if $n = m$, or zero otherwise). Integer translates of the sinc possess an orthonormal property:[11]

$$n \in \mathbb{Z} \quad \Rightarrow \quad \int_{-\infty}^{+\infty} \text{sinc}(z)\,\text{sinc}(z-n)\,dz = \text{sinc}(n) = \delta_{n0}. \tag{20}$$

Extending $\rho_i$ as $\hat{\rho}_i(z) = \rho_i(z)$ if $0 < z \leq \left|\partial^O \Pi(\mathbf{x}_i)\right|$ or $\hat{\rho}_i(z) = 0$ otherwise,

$$\int_0^{\left|\partial^O \Pi(\mathbf{x}_i)\right|} H_{ij}^O(z)\,\rho_i^O(z)\,dz = \int_{-\infty}^{+\infty} \hat{\rho}_i^O(z)\,\text{sinc}\left(\frac{z-z_j}{\Delta z}\right)dz \qquad \text{under } \textbf{(H)}, \tag{21}$$

with $z = 0$ at a corner of $\partial^O \Pi(\mathbf{x}_i)$. Approximating $\hat{\rho}_i^O(z)$ by cardinal interpolation,

$$\hat{\rho}_i^O(z) = \sum_{k \in \mathcal{S}_i^O} \rho_i^O(z_k)\,H_{ik}^O(z) + \varepsilon(\Delta z), \tag{22}$$

with $\varepsilon(\Delta z) \to 0$ when $\Delta z \to 0^+$. The off-diagonal, nonzero entries $G_{ij}$ are the sum of integrals like (21), for which

$$\lim_{\Delta z \to 0^+} \int_{-\infty}^{+\infty} \hat{\rho}_i^O(z)\text{sinc}\left(\frac{z-z_j}{\Delta z}\right)dz = \sum_{k \in \mathcal{S}_i^O} \int_{-\infty}^{+\infty} \rho_i^O(z_k)\,\text{sinc}\left(\frac{z-z_k}{\Delta z}\right)\text{sinc}\left(\frac{z-z_j}{\Delta z}\right)dz$$

$$= \sum_{k \in \mathcal{S}_i^O} \Delta z \rho_i^O(z_k) \int_{-\infty}^{+\infty} \text{sinc}(\xi)\,\text{sinc}\left[\xi - \left(\frac{z_j}{\Delta z} - \frac{z_k}{\Delta z}\right)\right]d\xi$$

$$= \sum_{k \in \mathcal{S}_i^O} \Delta z \rho_i^O(z_k)\,\text{sinc}\left(\frac{z_j}{\Delta z} - \frac{z_k}{\Delta z}\right)$$

$$= \Delta z \rho_i^O(z_j) \geq 0, \tag{23}$$

because $\rho_i^O(z)$ is nonnegative—and zero outside of the patch interface—and by (20). Then, by (16), one such $G_{ij}$ is a sum of nonpositive numbers. $\square$

---

[11] We have not found a reference for this well known fact and give the proof here. Let $n \in \mathbb{Z}$ and $\mathcal{F}[f(z)](\omega) = \int_{-\infty}^{+\infty} f(z)e^{-2\pi i w z}\,dz$. Then, the l.h.s. of (20) equals $\mathcal{F}[\text{sinc}(z)\,\text{sinc}(z-n)](\omega) := F(\omega)$ evaluated at $\omega = 0$. By the convolution (denoted by *) and shift properties of Fourier transforms, $F(\omega) = \mathcal{F}[\text{sinc}(z)](\omega) * e^{2\pi i \omega n}\mathcal{F}[\text{sinc}(z)](\omega)$. Since $\mathcal{F}[\text{sinc}(z)] = \text{rect}(\omega)$ (where $\text{rect}(\omega) = 1$ if $|\omega| > \frac{1}{2}$, $\text{rect}(\omega) = \frac{1}{2}$ if $|\omega| = 1/2$, and $\text{rect}(\omega) = 0$ if $|\omega| < \frac{1}{2}$), one has $F(0) = -\int_{-1/2}^{+1/2} e^{-2\pi i n \xi}\,d\xi = \text{sinc}(n) = \delta_{n0}$.

**Theorem 6** *Neglecting interpolation errors, $G^*$ and $G$ are diagonally dominant by rows, with at least one row strictly so.*

*Proof.* According to the splitting (17), it suffices to prove the statement for $G^*$:

$$\sum_{1 \le j \ne i \le N-m} |G^*_{ij}| \le \sum_{1 \le j \ne i \le N} |G_{ij}| = \sum_{1 \le j \ne i \le N} G_{ij}$$

$$= \sum_{O \in \{E,N,W,S\}} \int_0^{\left|\partial^O \Pi(\mathbf{x}_i)\right|} \rho_i^O(z) \left( \sum_{1 \le j \ne i \le N} H_{ij}^O(z) \right) dz \le 1 = |G^*_{ii}| \quad (24)$$

by Theorem 5, since $\sum_{1 \le j \ne i \le N} H_{ij}^O(z) = \sum_{j \in \mathcal{S}_i^O} H_{ij}^O(z) \to 1$ when interpolation errors are dropped, by (14). By (15), the strict inequality takes place at least on any row $i$ such that $\Pi(\mathbf{x}_i)$ is a perimeter patch exposed to Dirichlet BCs, which must exist since the elliptic BVP is well posed. □

For convenience, the following fact is adapted from Theorem 8.2.27 in [9].

**Lemma 7** *Let $A = [a_{ij}]$ be a real, irreducible, diagonally dominant square matrix of order $N$. If there is $1 \le i \le N$ such that $|a_{ii}| > \sum_{1 \le j \ne i \le N} |a_{ij}|$, then $A$ is invertible.*

**Corollary 8** *$G$ is an invertible M-matrix.*

*Proof.* By Theorem 6 and Lemma 7, $G^*$ is invertible, whence $G$ is invertible by Lemma 1. Since $G = I - B$ and $B$ is nonnegative by Theorem 5, one must only prove that the spectral radius of $B$ is less than one[12]. By Gershgorin's circles theorem it holds that $\rho(B) \le 1$, because the diagonal elements of $B$ are zero and $\sum_{i \ne j=1}^N |B_{ij}| = \sum_{i \ne j=1}^N B_{ij} \le 1$ by Theorem 6. On the other hand, since $B$ is nonnegative, the Perron-Frobenius theorem dictates that $\rho(B)$ is the largest eigenvalue of $B$, associated to the Perron eigenvector, $\vec{p}$, of $B$. Since $G\vec{p} = (1 - \rho(B))\vec{p}$, and $G$ is invertible, $\rho(B) \ne 1$. Therefore, $\rho(B) < 1$. □

### 4.3 Bounding the norm of the inverse

First, we check that nonnegative matrices preserve componentwise inequalities.

**Lemma 9** *Let $A \ge 0$ be an $N \times N$ matrix and $\vec{x} \ge \vec{y}$ componentwise. Then $A\vec{x} \ge A\vec{y}$.*

*Proof.* Since $A_{ij} \ge 0$ and $x_k \ge y_k$ ($1 \le i, j, k \le N$), it holds that $A_{ij}x_k \ge A_{ij}y_k$. Then

$$(A\vec{x})_i = \underbrace{A_{i1}x_1}_{\ge A_{i1}y_1} + \underbrace{A_{i2}x_2}_{\ge A_{i2}y_2} + \ldots + \underbrace{A_{iN}x_N}_{\ge A_{iN}y_N} \ge (A\vec{y})_i.$$ □

---
[12]https://nhigham.com/2021/03/16/what-is-an-m-matrix/

**Theorem 10** *Consider the elliptic BVP* (1) *with possibly both absorbing and reflecting BCs, and a PDDSparse discretisation into N knots, none of which lies at $\partial\Omega_A \cap \partial\Omega_R$ (the intersection would take place in the middle of a perimeter patch). Then,*

$$\|G^{-1}\|_\infty \leq 1 + \frac{\max_{\mathbf{x}_i \in \Omega} \mathbb{E}[\tau_i^\Omega]}{\min_{\mathbf{x}_j \in (\overline{\Omega}\setminus\partial\Omega_A)} \mathbb{E}[\tau_j^\#]}. \tag{25}$$

*Proof.* For a BVP like (1), $G$ only depends on $\mathcal{L}$, $c$, and $\varphi$, as well as the PDDSparse discretisation—not on $f, g, \psi$. Consider the two Poisson BVPs:

$$\begin{cases} \mathcal{L}w_s(\mathbf{x}) + \left(\min_{\mathbf{x}_j \in (\overline{\Omega}\setminus\partial\Omega_A)} \mathbb{E}[\tau_j^\#]\right)^{-1} = 0 \text{ in } \Omega, \\ \mathbf{N}^\top \nabla w_s(\mathbf{x}) = 0 \text{ on } \partial\Omega_R, \quad w_s(\mathbf{x}) = s = \{0, 1\} \text{ on } \partial\Omega_A. \end{cases} \tag{26}$$

Note that both BVPs give rise to the same $G$, and $w_1(\mathbf{x}) = 1 + w_0(\mathbf{x})$ in $\overline{\Omega}$. By (3),

$$w_0(\mathbf{y}) = \mathbb{E}\left[\int_0^{\tau^\Omega(\mathbf{y})} \left(\min_{\mathbf{x}_j \in (\overline{\Omega}\setminus\partial\Omega_A)} \mathbb{E}[\tau_j^\#]\right)^{-1} dt\right] = \left(\min_{\mathbf{x}_j \in (\overline{\Omega}\setminus\partial\Omega_A)} \mathbb{E}[\tau_j^\#]\right)^{-1} \mathbb{E}[\tau^\Omega(\mathbf{y})] \tag{27}$$

and therefore, for $\mathbf{x} \in \overline{\Omega}$,

$$w_1(\mathbf{x}) \leq 1 + \frac{\max_{\mathbf{x}_i \in \Omega} \mathbb{E}[\tau_i^\Omega]}{\min_{\mathbf{x}_j \in (\overline{\Omega}\setminus\partial\Omega_A)} \mathbb{E}[\tau_j^\#]} \geq 2. \tag{28}$$

In the denominator above, the knots on the Dirichlet boundary have been excluded. Along with the hypothesis, this ensures that the minimum is strictly positive. Note also that $\mathbb{E}[\tau_i^\Omega] \geq \mathbb{E}[\tau_i^\#]$ for any $1 \leq i \leq N$.

Now, let the PDDSparse nodal solution of (26) with $s = 1$ be the vector $\vec{W}$, such that $(\vec{W})_k = w(\mathbf{x}_k)$ for $k = 1, \ldots, N$. It solves the linear system $G\vec{W} = \vec{q}$, with

$$(\vec{q})_k = \begin{cases} \mathbb{E}[Y_{\tau_k^\#} \mathbb{1}_{\{\mathbf{X}_{\tau_k^\#} \in \partial\Omega_A\}} | \mathbf{X}_0 = \mathbf{x}_k] + \left(\mathbb{E}[\tau_k^\#] / \min_{\mathbf{x}_j \in (\overline{\Omega}\setminus\partial\Omega_A)} \mathbb{E}[\tau_j^\#]\right) \geq 1 & \text{if } \mathbf{x}_k \in \Omega \cup \partial\Omega_R, \\ 1 & \text{if } \mathbf{x}_k \in \partial\Omega_A. \end{cases}$$

Then $G\vec{W} \geq \vec{1}$, where $\vec{1}$ is a vector of $N$ ones, and, since $G^{-1} \geq 0$,

$$\vec{W} \geq G^{-1}\vec{1}, \tag{29}$$

where Lemma 9 has been invoked. By elementary calculations, (28) and (29),

$$\|G^{-1}\|_\infty = \max_{1 \leq i \leq N} \left(|G^{-1}|\vec{1}\right)_i = \max_{1 \leq i \leq N} \left(G^{-1}\vec{1}\right)_i \leq \|\vec{W}\|_\infty$$

$$= \max_{1 \leq i \leq N} |(\vec{W})_i| \leq \max_{\mathbf{x} \in \overline{\Omega}} w_1(\mathbf{x})$$

$$\leq 1 + \left(\max_{\mathbf{x}_i \in \Omega} \mathbb{E}[\tau_i^\Omega] / \min_{\mathbf{x}_j \in (\overline{\Omega}/\partial\Omega_A)} \mathbb{E}[\tau_j^\#]\right). \qquad \square$$

## 4.4 The condition number of $G$

**Lemma 11** *Let H be the interface length. The condition number of G is bounded as*

i) *In the infinity norm:* $\kappa_\infty(G) \leq 2\|G^{-1}\|_\infty$.

ii) *In the spectral norm:* $\kappa_2(G) \leq O(1)\sqrt{nHN}\kappa_\infty(G) \leq N\kappa_\infty(G)$ *(in fact, $nH \ll N$).*

iii) *Skeel's condition number:* $\text{cond}_\infty(G) := \left\| |G^{-1}| \cdot |G| \right\|_\infty \leq 1 + \kappa_\infty(G)$.

*Proof.* The condition number in any given matrix norm $\|\cdot\|$ is $\kappa_{\|\cdot\|}(G) = \|G^{-1}\| \cdot \|G\|$. For *i)*, $\|G^{-1}\|_\infty$ is bounded by Theorem 10; what remains is:

$$\|G\|_\infty = \max_{1 \leq i \leq N} \sum_{j=1}^{N} |G_{ij}| = 1 + \max_{1 \leq i \leq N} \sum_{j \neq i} |G_{ij}| \leq 2. \tag{30}$$

Regarding *ii)*, the upper bound stems from applying to $G$ (and $G^{-1}$) the standard inequality $\|A\|_2 \leq \sqrt{N}\|A\|_\infty$. However, in this case it is better to use the (also standard) relation $\|G\|_2^2 \leq \|G\|_1 \cdot \|G\|_\infty$. Note that the nonempty entries along column $j^{th}$ of $G$ are $G_{jj} = 1$, and as many more as there are knots which have $\mathbf{x}_j$ in their stencils. Even minding elongation, they all sit on a few interfaces around $\mathbf{x}_j$. If the latter have length $H$, there are $O(1) \times nH$ such knots $\{\mathbf{x}_k\}$, for which $-1 \leq G_{kj} \leq 0$. Consequently, $\sum_{1 \leq i \leq N} |G_{ij}| \leq 1 + nHO(1)$, and this is a conservative bound for $\|G\|_1$. In sum, $\|G\|_2 \leq \sqrt{2}O(1)\sqrt{nH}$, and the sharper bound in *ii)* follows. Finally, *iii)* is a consequence of the fact that $G = 2I - |G|$. □

By Theorem 10, studying the condition number of $G$ hinges upon estimating mean FETs, which is problem dependent. We give next two preparatory lemmas for the Brownian motion. (The first one is well known, see [13, example 7.4.2].)

**Lemma 12** *The expected FET of the Wiener process from inside a circle of radius R, starting at a distance r from its centre at the origin, is $\mathbb{E}[\tau] = (R^2 - r^2)/2 =: E^O(R, r)$.*

**Lemma 13** *Assume $\mathcal{L}$ is the Laplacian, possibly multiplied by a constant. Then $\mathbb{E}[\tau^\Omega(\mathbf{x})]$ scales with the square of the diameter of the domain. Concretely,*

$$\Omega' := \{\mathbf{x}' := \Lambda \mathbf{x} : \mathbf{x} \in \Omega\} \quad \Rightarrow \quad \mathbb{E}[\tau^{\Omega'}(\mathbf{z})] = \Lambda^2 \mathbb{E}[\tau^\Omega(\mathbf{z}/\Lambda)]. \tag{31}$$

*Proof.* Suppose that $\Omega$ is stretched by a factor $\Lambda > 0$ into $\Omega'$. By (5), the mean FET from point $\mathbf{x}' \in \Omega'$ obeys the BVP

$$\Delta'\mathbb{E}[\tau^{\Omega'}(\mathbf{x}')] + 2 = 0 \text{ in } \Omega', \ \mathbb{E}[\tau^{\Omega'}(\mathbf{x}')] = 0 \text{ on } \partial\Omega'_A, \ \mathbf{N}^\top\nabla\mathbb{E}[\tau^{\Omega'}(\mathbf{x}')] = 0 \text{ on } \partial\Omega'_R, \tag{32}$$

where $\Delta'$ is the Laplacian in the stretched coordinates. Note that (5) can be rewritten (by just renaming variables) as

$$\Delta\mathbb{E}[\tau^\Omega(\mathbf{x}'/\Lambda)] + 2 = 0 \text{ in } \Omega', \ \mathbb{E}[\tau^\Omega(\mathbf{x}'/\Lambda)] = 0 \text{ on } \partial\Omega'_A, \ \mathbf{N}^\top\nabla\mathbb{E}[\tau^{\Omega'}(\mathbf{x}'/\Lambda)] = 0 \text{ on } \partial\Omega'_R.$$

On the other hand, by the chain rule, and regrouping constants

$$\Delta' \left( \Lambda^2 \mathbb{E}[\tau^{\Omega}(\mathbf{x}'/\Lambda)] \right) + 2 = 0 \text{ in } \Omega',$$
$$\left( \Lambda^2 \mathbb{E}[\tau^{\Omega}(\mathbf{x}'/\Lambda)] \right) = 0 \text{ on } \partial\Omega'_A, \quad \mathbf{N}^\top \nabla \left( \Lambda^2 \mathbb{E}[\tau^{\Omega'}(\mathbf{x}'/\Lambda)] \right) = 0 \text{ on } \partial\Omega'_R. \quad (33)$$

By inspection, the unique solution to (32) is that in (33). □

**Condition number rates of typical scenarios.** In the upcoming analysis, let $\Omega$ be a square of side $L$ ($|\Omega| = L^2$) partitioned into $M = m \times m$ square subdomains of side $H = L/m$. The interfaces are discretised according to a linear density of knots $n$. In order to draw mean FET estimates, we assume that $\mathcal{L}$ in (1) is half the Laplacian—so that the stochastic trajectories are Brownian motions. In that case, the minimiser in the denominator of (25) takes place on any knot sitting at a distance equal to the internodal separation, $\Delta z = 1/n$, of its patch boundary. One such closest knot, labeled as $\mathbf{x}_{min}$, is depicted in Figure 3 (right). From $\mathbf{x}_{min}$, the mean FET from inside its patch and from inside the inscribed circle are $\mathbb{E}[\tau^\#(\mathbf{x}_{min})]$ and $E^O(H/2, H/2 - \Delta z)$, respectively. Clearly, $\mathbb{E}[\tau^\#(\mathbf{x}_{min})] > E^O(H/2, H/2 - \Delta z)$. By Lemma 12, since $(H/2)^2 - (H/2 - \Delta z)^2 \approx H\Delta z$, one has $(\mathbb{E}[\tau^\#(\mathbf{x}_{min})])^{-1} \lesssim 2/(H\Delta z)$. Then, by Theorem 10 and Lemma 11, item i),

$$\kappa_\infty(G) \lesssim 2 + \frac{2 \max_{\mathbf{x}_i \in \Omega} \mathbb{E}[\tau_i^\Omega]}{H\Delta z}. \quad (34)$$

**Remark 3** *A sharper bound for $\mathbb{E}[\tau^\#(\mathbf{x}_{min})]$ can be obtained by solving (32) exactly by separation of variables in the patch. Nonetheless, the leading orders are the same. For $(x, y) \in [-H/2, H/2]^2$, the mean FET is $\mathbb{E}[\tau^\#(x, y)] = \frac{32 H^2}{\pi^4} \sum_{k,l=1}^{\infty} \frac{(-1)^{k+l} \cos\left(\frac{2k+1}{H}\pi x\right) \cos\left(\frac{2l+1}{H}\pi y\right)}{(2k+1)(2l+1)[(2k+1)^2 + (2l+1)^2]}.$*

**Remark 4** *The bound (34) is strikingly reminiscent of $\kappa_2(S) < C/(hH)$, valid for the Schur complement $S$ of FEM discretisations of diameter $h$ of Poisson equations with Dirichlet BCs [17, lemma 4.11]. There, $C$ is a constant independent of the discretisation, and indentifying $h$ with $\Delta z$ essentially yields (34), even though both domain decomposition methods and arguments are unrelated. While the norms are different, our experiments show that $\kappa_2(G)$ and $\kappa_\infty(G)$ are rather similar in practice.*

In addition to Remark 4, we shall bound the growth of $\kappa_\infty(G)$ when the total number of knots, $N$, increases according to three specific patterns (or *scenarios*):

i) *Larger domain, same accuracy.* The numerator of (34) grows proportional to $|\Omega| = L^2$, thus proportional to $N$. The denominator, on the other hand, is constant, as floating patches and knots therein remain the same. This is a typical *weak scalability scenario* in parallel computing, where the problem size grows in order to make a more efficient use of the available cores.

$$L \text{ grows while } H, n \text{ constant} \quad \Rightarrow \quad \kappa_\infty(G) \leq \mathcal{O}(N). \quad (35)$$

ii) *Same domain, better accuracy.* Another weak scalability scenario where $\Omega$ and $H$ are fixed, while $n$ is refined seeking better numerical accuracy. Consequently, the numerator of (34) remains constant, but the denominator decreases as $\Delta z = 1/n$. In this scenario, $N \propto n$, and therefore

$$n \text{ grows while } L, H \text{ constant} \quad \Rightarrow \quad \kappa_\infty(G) \leq O(N). \quad (36)$$

iii) *Same domain and accuracy, more subdomains.* Here $L$ and $n$ are constant, but $H$ shrinks. Even though $\Delta z$ is constant too, the denominator of (34) decreases because the patches are smaller. By (34) and the fact that—in this scenario—$N \approx 2mLn = 2L^2/(H\Delta z)$, we have $(H\Delta z)^{-1} \propto N$, and again

$$m \text{ grows } (H \text{ shrinks}), L, n \text{ constant} \quad \Rightarrow \quad \kappa_\infty(G) \leq O(N). \quad (37)$$

The fact that $\kappa_\infty(G)$ grows moderately in typical discretisation setups is a measure of stability. $G$ can also be very efficiently preconditioned, as shown next.

## 5 A preconditioner for PDDSparse

In an actual PDDSparse simulation, the expected values (11)-(12) making up the entries of $G$ and $\vec{b}$ are produced by averaging a large number of SDE realisations, $\tilde{N}$, approximated for instance with an Euler-Maruyama scheme with timestep $h$. This incurs a bias linked to $h$, plus the (random) statistical error of replacing the expectations by means. (See [3] for details.) Additionally, the finite discretisation preventing the interpolation cardinal functions from being actual sincs, may make $\tilde{G}$ depart slightly from a diagonally dominant M-matrix. However, since accuracy dictates that $\tilde{N} \gg 1$ and $0 < h \ll 1$, every realised matrix can be construed as a random perturbation of the "underlying", deterministic, M-matrix $G$. In other words, $\tilde{G}(\omega) \approx G$ can be expected.

### 5.1 Naïve Neumann preconditioner

Given a realisation of $\tilde{G}$, we seek a perturbation matrix $\tilde{E}$ such that

$$\tilde{G} + \tilde{E} := P \text{ where } \begin{cases} P \text{ is a diagonally dominant M-matrix, and} \\ \|\tilde{E}\|/\|\tilde{G}\| \ll 1 \text{ in a suitable matrix norm.} \end{cases}$$

We can easily tailor a satisfactory matrix $\tilde{E}$ to $\tilde{G}$ as follows.

**Algorithm (A).** Given $\tilde{G}$, let $\tilde{E}_{ij} = 0$ for $1 \leq i, j \leq N$, and then:

I) Get the off-diagonal positive terms of $\tilde{G}$: $\tilde{E}_{ij} = -\tilde{G}_{ij}$ iff $i \neq j$ and $\tilde{G}_{ij} > 0$.

II) Enforce diagonal dominance of $\tilde{G} + \tilde{E}$: $\tilde{E}_{ii} = \delta$ for $1 \leq i \leq N$, where

$$\delta = \max_{1 \leq i \leq N} \left\{ 0, \left( \sum_{j=1}^{N} \tilde{G}_{ij} + \tilde{E}_{ij} \right) - 1 \right\}. \quad (38)$$

III) (Optional) If $\delta$ is unsatisfactorily large, let $i = \alpha$ be the row where the max above takes place. Recompute the entries $G_{\alpha j}$, and repeat steps I) and II).

**Remark 5** *As we shall see, the preconditioning of $\tilde{G}$ is computationally enhanced by having as small a $\delta$ as possible. The max value $\delta$ returned by Algorithm (A) after steps I) and II) is a random variable. It can occur that the max is an outlier with respect to the histogram of $\delta_i := -1 + \sum_{j=1}^{N}(\tilde{G}_{ij} - \tilde{E}_{ij})$, which can be generated on the fly. In that case, it might be worthwhile to recompute only that row's entries from scratch (this is a privilege of PDDSparse). The decision depends on balancing the cost of recomputing a row against the expected improvement. To estimate the latter, one can resort to basic statistics on the histogram of $\{\delta_i\}$, such as Markov's or Chebyshev's inequalities.*

According to Algorithm (A), the *Neumann preconditioner P* is

$$P = \tilde{G} + \tilde{E}_\delta =: (1 + \delta)I - \tilde{B}_\delta, \tag{39}$$

(this new notation highlights that $\tilde{E}_\delta$ and $\tilde{B}_\delta$ depend[13] on $\delta$) where $(\tilde{B}_\delta)_{ii} = 0$ and $(\tilde{B}_\delta)_{i \neq j} \geq 0$. By mimicking the proof of Corollary 8, it holds that $\rho(\tilde{B}_\delta) < 1 + \delta$, so that $P$ is indeed an M-matrix. Then, the Neumann series of $P^{-1}$ converges:

$$P^{-1} = \frac{1}{1+\delta}\left(I - \frac{\tilde{B}_\delta}{1+\delta}\right)^{-1} = \frac{1}{1+\delta}\left(I + \sum_{k=1}^{\infty}\frac{(\tilde{B}_\delta)^k}{(1+\delta)^k}\right) = \frac{1}{1+\delta}\sum_{k=0}^{\infty}\frac{(\tilde{B}_\delta)^k}{(1+\delta)^k}. \tag{40}$$

In any submultiplicative matrix norm such that $\|I\| = 1$, it holds

$$\|P^{-1}\tilde{G}\| = \|P^{-1}(P - \tilde{E}_\delta)\| \leq 1 + \|P^{-1}\| \cdot \|\tilde{E}_\delta\|, \tag{41}$$

$$\|(P^{-1}\tilde{G})^{-1}\| = \|\tilde{G}^{-1}(\tilde{G} + \tilde{E}_\delta)\| \leq 1 + \|\tilde{G}^{-1}\| \cdot \|\tilde{E}_\delta\|. \tag{42}$$

Consequently, the corresponding preconditioned condition number obeys

$$\kappa_{\|\cdot\|}(P^{-1}\tilde{G}) \leq 1 + \left(\|\tilde{G}^{-1}\| + \|P^{-1}\|\right)\|\tilde{E}_\delta\| + \left(\|\tilde{G}^{-1}\| \cdot \|P^{-1}\|\right)\|\tilde{E}_\delta\|^2. \tag{43}$$

Above, $\|P^{-1}\| \approx \|G^{-1}\|$ is controlled by Theorem 10, and $\|\tilde{G}^{-1}\|$ is finite—preconditioning is futile otherwise. In the infinity norm, $\|\tilde{E}_\delta\|_\infty = \mathcal{O}(\delta)$ and it is expected that $\delta \ll 1$, implying that $\kappa_\infty(P^{-1}\tilde{G}) = 1 + \mathcal{O}(\delta)$.

Of course, the Neumann series (40) must be truncated in practice. Let

$$P_t^{-1} := \frac{1}{1+\delta}\sum_{k=0}^{t}\frac{(B_\delta)^k}{(1+\delta)^k}. \tag{44}$$

Since

$$P^{-1} - P_t^{-1} = \frac{1}{1+\delta}\sum_{k=t+1}^{\infty}\frac{(\tilde{B}_\delta)^k}{(1+\delta)^k} = \frac{\tilde{B}_\delta^{t+1}}{(1+\delta)^{t+1}}P^{-1}, \tag{45}$$

---

[13]On the other hand, we refrain from writing $P_\delta$ for later notational simplicity.

and because $P^{-1}$ commutes with $\tilde{B}_\delta$ (hence with every power), we have

$$P_t^{-1} = \left(I - \frac{\tilde{B}_\delta^{t+1}}{(1+\delta)^{t+1}}\right)^{-1} P^{-1} = P^{-1}\left(I - \frac{\tilde{B}_\delta^{t+1}}{(1+\delta)^{t+1}}\right)^{-1}. \qquad (46)$$

Letting $T_{t+1} := I - \tilde{B}_\delta^{t+1}/(1+\delta)^{t+1}$ and replicating formulas (41)-(42), the truncated preconditioned condition number obeys

$$\kappa_{\|\cdot\|}(P_t^{-1}\tilde{G}) \leq \|T_{t+1}\| \cdot \|T_{t+1}^{-1}\| \left[1 + \left(\|\tilde{G}^{-1}\| + \|P_t^{-1}\|\right)\|\tilde{E}_\delta\| + O\left(\|\tilde{E}_\delta\|^2\right)\right], \qquad (47)$$

which goes to $\|T_{t+1}\| \cdot \|T_{t+1}^{-1}\| = \kappa_{\|\cdot\|}\left(I - \tilde{B}^{t+1}/(1+\delta)^{t+1}\right)$ as $\|\tilde{E}_\delta\| \to 0^+$.

Equation (47) provides guidance on picking $t$. In order to estimate $\|T_{t+1}^{-1}\|$, one can use the Neumann series again, since $\rho\left([\tilde{B}_\delta/(1+\delta)]^{t+1}\right) < \rho\left(\tilde{B}_\delta/(1+\delta)\right) < 1$. Then, if $t$ is such that $(1+\delta)^{-(t+1)} \gg (1+\delta)^{-2(t+1)}$, and cutting at $k = 1$, one has

$$T_{t+1}^{-1} = I + \sum_{k=1}^{\infty}\left(\frac{\tilde{B}_\delta^{t+1}}{(1+\delta)^{t+1}}\right)^k \approx 2I - T_{t+1}, \qquad (48)$$

whence estimates such as $\kappa_{\|\cdot\|}(T_{t+1}) \lesssim 2\|T_{t+1}\|(1 + \|T_{t+1}\|)$ are straightforward.

Regarding implementation, the truncated Neumann preconditioner (44) is very well suited to Krylov methods for solving the actual system $\tilde{G}\vec{u} = \vec{b}$, for any matrix-vector product $\vec{z} = P_t^{-1}\vec{x}_0$ is the result of $t$ iterations of the recursion

$$\vec{x}_k = \vec{x}_0 + (1+\delta)^{-1}\tilde{B}_\delta \vec{x}_{k-1}, \qquad (49)$$

followed by $\vec{z} = (1+\delta)^{-1}\vec{x}_t$, where $\tilde{B}_\delta\vec{x}_{k-1}$ involves only $O(nN)$ operations, since $\tilde{B}_\delta$ is sparse with bandwidth $O(n)$. (In realistic applications, $N \gg n$.)

## 5.2 A much better Neumann-Arnoldi preconditioner

Experimentally, $P_t^{-1}$ turns out to be effective in reducing both the condition number and the number of GMRES [14] iterations for solving $P_t^{-1}\tilde{G}\vec{u} = P_t^{-1}\vec{b}$. Unfortunately though, the *total computational time* is only marginally better or—more often—distinctly worse.

Since $G \approx \tilde{G}$ is unsymmetric, there is no *a priori* link between the condition number and the speed of convergence of GMRES. Instead, the key for fast convergence is that the spectrum of the preconditioned matrix be clusterised around the point $(1,0)$ on the complex plane [15].

In order to enhance the spectrum of $P_t^{-1}\tilde{G}$ while keeping $t$ low, we have borrowed an idea from the recent paper [18] by Zheng, Xi, and Saad[15]. Let

$$\tilde{G}^{-1} = P_t^{-1} + R, \qquad (50)$$

---

[14]Recall that $G \approx \tilde{G}$ is not symmetric, hence ruling out conjugate gradient.
[15]In [18], $C_0$ is motivated by a Schur complement matrix with overlapping block diagonal structure, unrelated to M-matrices.

where $R$ can be thought of as the remainder. Then,

$$I = \tilde{G}P_t^{-1} + \tilde{G}R \tag{51}$$

and

$$\tilde{G}^{-1} = P_t^{-1}(I - \tilde{G}R)^{-1}. \tag{52}$$

The insight of [18] is to approximate the rightmost inverse by a low-rank proxy. Let $0 < r \ll N$ and $I_r$ be the identity in $\mathbb{R}^{r \times r}$. First, construct $H_r$ and $V_r$ in

$$\tilde{G}R = I - \tilde{G}P_t^{-1} \approx V_r H_r V_r^\top, \tag{53}$$

where $H_r$ is an $r \times r$ Hessenberg matrix and the $N \times r$ matrix $V_r$ has orthonormal columns. Both are simultaneously generated after $r$ steps of the incomplete Arnoldi procedure started with some vector $\vec{y}_0 \neq \vec{0}$. Importantly, this can be carried out in terms of matrix-vector products only, via (49), so that is suitable for a distributed computer. Then, using the Sherman-Morrison formula gives

$$(I - \tilde{G}R)^{-1} \approx (I - V_r H_r V_r^\top)^{-1} = I + V_r \left[(I_r - H_r)^{-1} - I_r\right] V_r^\top. \tag{54}$$

The point is that $(I_r - H_r)^{-1}$ is small enough that it can be inverted explicitly and stored. Plugging (54) into (52) defines the *Neumann-Arnoldi preconditioner*:

$$\Pi_{t,r} = P_t^{-1}\left(I + V_r \left[(I_r - H_r)^{-1} - I_r\right] V_r^\top\right). \tag{55}$$

As will be seen in Section 6, the spectrum of $\Pi_{t,r}\tilde{G}$ fits in a much smaller circle around $1 + 0i$ than that of $\tilde{G}$ (and of $P_t^{-1}\tilde{G}$), resulting in significantly fewer GMRES iterations. Let one matrix-vector product $\tilde{B}_\delta \vec{x}$ (or $\tilde{G}\vec{x}$, with similar sparsity) be the unit of cost, and #$it_0$ and #$it_{t,r}$ the number of GMRES iterations to convergence. By (49), the action of $P_t^{-1}$ costs $t$. Then, a preconditioned GMRES solve essentially costs #$it_{t,r} \times t$ (the products in (49) per iteration), plus the initial cost of generating $V_r$ and $H_r$ in (53). Disregarding the $O(r^3)$ cost of the matrix inversion (54), the bulk of the overhead is $r$ times the matrix-vectors in the incomplete Arnoldi procedure for (53), i.e. $rt$. Consequently, the Neumann-Arnoldi preconditioner pays when $t(r + $#$it_{t,r}) < $#$it_0$. For instance, this would occur if $t = O(1)$ and the iterations drop #$it_0 - $#$it_{t,r}$ is large compared to $r$.

## 6  Numerical experiments

We solved Poisson's equation $\nabla^2 u = F$, with Dirichlet BCs and exact solution

$$u(x,y) = 3 + \tfrac{1}{3}\sin\left(\sqrt{1 + \tfrac{x^2}{100} + \tfrac{y^2}{50}}\right) + \tfrac{1}{3}\tanh\left[\sin\left(\tfrac{3x}{25} + \tfrac{y}{20}\right) + \sin\left(\tfrac{x}{20} - \tfrac{3y}{25}\right)\right] \tag{56}$$

(same as BC), on several square domains $\Omega = [-L/2, L/2]^2$. The code was written in C++/ OpenMP/ MPI, and run on the Finnish supercomputer LUMI in order to generate matrices $\tilde{G}(\omega)$. They were later analysed using MATLAB.

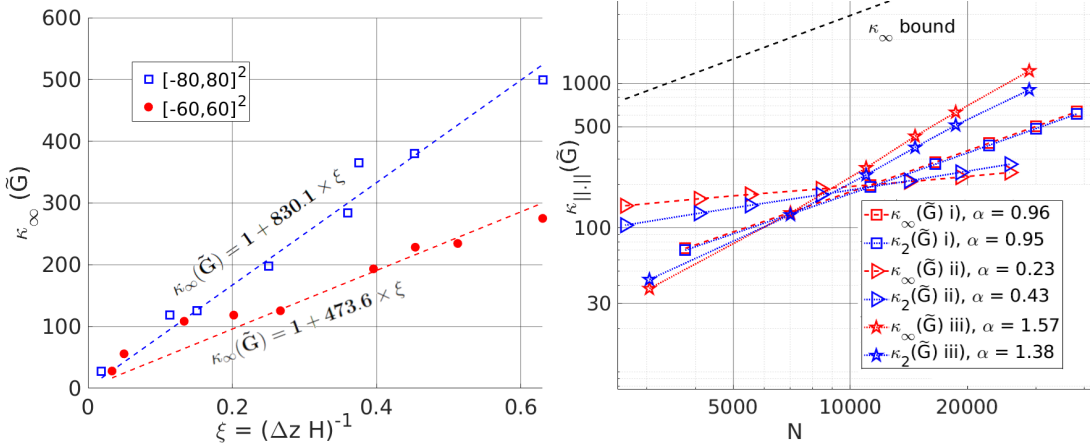

Figure 5: Numerical check of (34) (left) and of scenarios (35)-(37) (right, in log-log scale). Beware the different scales and the fact that $\tilde{G}$, not $G$, is used. In all scenarios i)-iii), the bound (34) is $\kappa_\infty(\tilde{G}) \lesssim 2 + 0.295 \times N$ (black dashed line).

**Raw (unpreconditioned) condition number.** Figure 5 (left) depicts the observed[16] $\kappa_\infty(\tilde{G})$ vs. $\xi := \Delta z H$ for $L = \{120, 160\}$, in good agreement with (34). Note that the fitted slopes are proportional to $L^2 = |\Omega| \propto \max_{\mathbf{x}\in\Omega} \mathbb{E}[\tau^\Omega]$. Results on the right pertain to the three scenarios described in (35)-(37). The fitted exponents of $\kappa_\infty(\tilde{G}) - 1 = \beta N^\alpha$ are reported. Only in scenario iii) does $\kappa_\infty(\tilde{G})$ grow superlinearly with $N$—but well below the theoretical upper bound. Intriguingly, it turns out that $\kappa_\infty(\tilde{G}) \approx \kappa_2(\tilde{G})$, but we cannot explain why.

**Neumann preconditioner.** Next, we consider the system $P_t^{-1}\tilde{G}\vec{u} = P_t^{-1}\vec{b}$, where $P_t^{-1}$ is given by (44) and implemented as (49). As shown on the left of Figure 6, the Neumann preconditioner is effective even when $\delta$ is relatively large (which corresponds to largish statistical errors). Observations are aligned with the theory; in particular the bound (47), depicted by the segments on each data point. While the number of GMRES iterations does decrease with $t$ (check Table 1), the effect diminishes rapidly with $t$, too. This is due to the relatively poor clusterisation of eigenvalues achieved by $P_t^{-1}$ (see Figure 6, top right).

**Neumann-Arnoldi preconditioner.** The low-rank-enriched preconditioner $\Pi_{t,r}$ in (55) dramatically improves the situation, in terms of both the preconditioned spectrum (Figure 6, bottom right), and of GMRES iterations[17].

Table 1 gives details on a battery of tests with our largest matrix ($N = 94580$), including the wallclock times of GMRES and of producing (once) the low-rank correction (53) and (54), on a desktop in MATLAB. We stress that those times—and hence the optimal parameters—are problem-, hardware-, and implementation-dependent. Regarding aggregate computational times, we could not get outstanding speedup in any case. The reason is as follows. According to the closing analysis in Section 5, and given that in our examples the optimal $r$ is about 100 (where $100 \ll N$), showing clear acceleration would take $\#it_0 - \#it_{r,t,1} \gg 100$, but this is impossible since $\#it_0 \approx 200$. In sum, even our biggest example was not nearly big enough to make a difference[18]. This is actually very encouraging, for it shows that the system will likely remain stable and preconditionable when truly huge simulations will be tackled.

---

[16] Recall that $G$ is a theoretical construct, and the theory is validated on the "noisy" matrix $\tilde{G}$.
[17] The preconditioned condition number actually grows with this preconditioner.

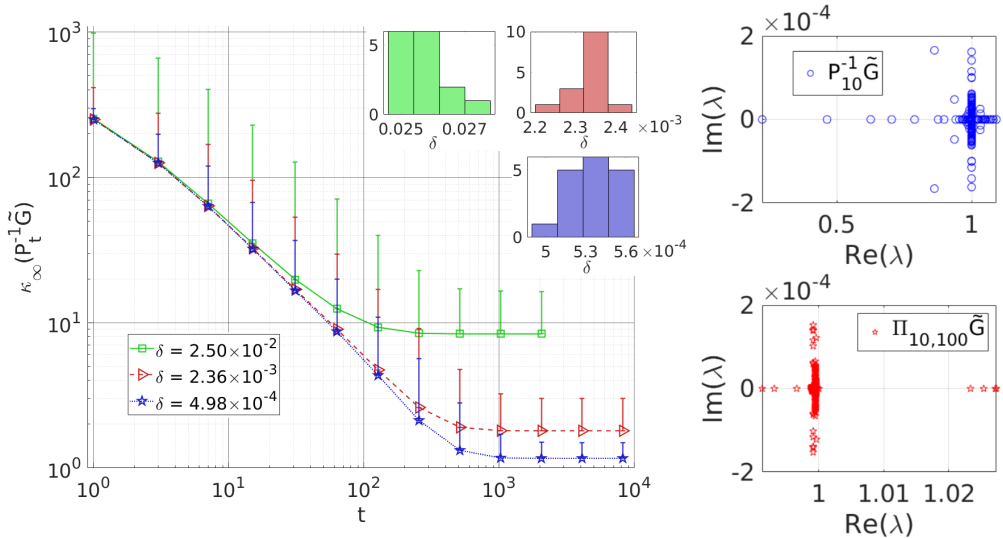

Figure 6: Left: effect of $t$ and $\delta$ on the preconditioned condition number of $P_t^{-1}\tilde{G}$. The representativity of the three exemplary $\delta$'s can be assessed by comparing with the histograms of the sample populations from which they are drawn (insets). Each histogram are several statistically independent matrices $\tilde{G}(\omega)$ generated with same values of $h, \tilde{N}$ of Euler-Maruyama. Top right: one full spectrum of $P_{10}^{-1}\tilde{G}(\omega = \omega_0)$. Bottom right: same matrix, spectrum for $\Pi_{10,100}^{-1}\tilde{G}(\omega = \omega_0)$. (Note the different scale of either spectrum.)

The preconditioned GMRES iterations are systematically fewer when the initial vectors $\vec{y}_0$ and $\vec{x}_0$—for the incomplete Arnoldi and GMRES algorithms respectively—are coupled as $\vec{y}_0 = \vec{b} - \tilde{G}\vec{x}_0$. This remarkable fact can be intuitively explained assuming that $\delta = 0$ and setting $t = 0$. Then, $\tilde{G} = G$, $P_t^{-1} = I$, and $GR = B$, of $G = I - B$. There are three Krylov spaces: $\mathcal{V}_r^A = \text{span}\{\vec{y}_0, B\vec{y}_0, \ldots, B^r\vec{y}_0\}$ for (incomplete) Arnoldi; $\mathcal{V}_K^{GMRES} = \text{span}\{\vec{r}_0, \Pi_{0r}G\vec{r}_0, \ldots, (\Pi_{0r}G)^K\vec{r}_0\}$ for GMRES after $K$ preconditioned iterations; and $\vec{u}$ lives in $\mathcal{V}_\infty^{ex} = \text{span}\{\vec{b}, B\vec{b}, B^2\vec{b} \ldots\}$ because of the M-matrixness of $G$. In this situation, $\vec{x}_0 = \tilde{G}^{-1}(\vec{y}_0 - \vec{b})$, $\vec{r}_0 = \Pi_{0r}(G\vec{x}_0 - \vec{b}) = -\Pi_{0r}\vec{y}_0$, and—since $\vec{y}_0 \in \mathcal{V}_r^A$ and $\Pi_{0r} = I + V_r[(I_r - H_r)^{-1} - I_r]V_t^\top$ by (55)—it holds that $\vec{r}_0 \in \mathcal{V}_r^A$, and, in general, $\Pi_{0r}G\vec{r}_0 = \Pi_{0r}(I - B)\vec{r}_0 \in \mathcal{V}_{r+1}^A$. In other words, the GMRES iterations simply extend the space where the GMRES solution lives from $\mathcal{V}_r^A$ to $\mathcal{V}_{r+K}^A$. Then, the ideal case is when $\mathcal{V}_{r+K}^A = \text{span}\{\vec{b}, B\vec{b}, \ldots, B^{r+K}\vec{b}\}$, which occurs by taking $\vec{y}_0 = \vec{b}$ and $\vec{x}_0 = \vec{0}$. (However, picking a random $\vec{x}_0$ takes just one more iteration as long as $\vec{y}_0 = \vec{b} - \tilde{G}\vec{x}_0$.)

---

[18] Unfortunately, producing the required matrices was far beyond our usage constraints of LUMI.

# 7 Conclusions

Despite lacking symmetry and its complex error structure, the PDDSparse matrix is underpinned by fruitful properties, unveiled by the analysis presented in this article. Thanks to them, we have carried out a theoretical study of the raw as well as preconditioned condition numbers, using a manufactured, yet simple and general, preconditioning scheme. The predicted and observed convergence rates of the raw and preconditioned condition numbers, as well as of GMRES iterations, are quite encouraging. Results presented here pertain not only to the seminal paper [3], but to the improved variant [11] as well.

In order to tackle truly massive simulations, the PDDSparse system must be itself substructured, as discussed in [3]. That work, currently under way, builds on this article—for the Schur complement preserves M-matrixness.

# Acknowledgements

We acknowledge grants 2018-T1/TIC-10914; 2022-5A/TIC-24233 of the Madrid Regional Government; PID2020-115088RB-I00 of Spanish AEI; and the EuroHPC Joint Undertaking for awarding us access to LUMI at CSC, Finland (ID: EHPC-DEV-2023D04-008).

Table 1: Raw ($t = r = 0$) and preconditioned GMRES iterations (with tolerance $10^{-12}$) for $\tilde{G}\vec{u} = \vec{b}$ ($t > 0, r = 0$ is the Neumann preconditioner). $\tilde{G}$ is of order $N = 94580$, with $\Omega = [-140, 140]^2$ into $28 \times 28$ subdomains. "Time LR" is that of Arnoldi and explicit inversion in (53) and (54). Figures in paretheses use $\vec{y}_0 = \vec{1}$ and $\vec{x}_0 = \vec{0}$; otherwise the starting vectors of Arnoldi and GMRES are coupled as $\vec{y}_0 = \vec{b}, \vec{x}_0 = \vec{0}$. In bold, the least total time we were able to achieve.

| $t$ | $r$ | #$it$ GMRES | time LR | time GMRES | total time |
|---|---|---|---|---|---|
| 0 | 0 | 193 | 0 | 13.21 | 13.21 |
| 0 | 10 | 179 (150) | 0.47 (0.43) | 11.88 (9.5) | 12.35 (9.93) |
| 0 | 25 | 161 (142) | 1.10 (1.12) | 10.42 (9.45) | 11.52 (10.67) |
| 0 | 50 | 134 (146) | 2.31 (2.30) | 8.35 (9.38) | 10.62 (11.68) |
| 0 | 75 | 109 (139) | 3.61 (3.58) | 6.81 (9.02) | 10.42 (12.60) |
| 0 | 100 | 84 (136) | 4.97 (4.89) | 5.11 (8.98) | 10.08 (13.87) |
| 1 | 0 | 105 | 0 | 9.33 | 9.33 |
| 1 | 10 | 80 | 0.76 (0.76) | 8.43 (7.00) | 9.19 (**7.76**) |
| 1 | 25 | 77 (85) | 1.93 (1.91) | 6.95 (7.53) | 8.88 (9.44) |
| 1 | 50 | 52 (77) | 3.99 (3.92) | 4.55 (6.89) | 8.54 (10.81) |
| 1 | 75 | 27 (67) | 6.05 (6.00) | 2.32 (6.03) | 8.37 (12.03) |
| **1** | **85** | **17** (61) | **6.81** (6.84) | **1.48** (5.50) | **8.29** (12.37) |
| 1 | 100 | 2 (61) | 8.29 (8.12) | 0.25 (5.50) | 8.54 (13.62) |
| 2 | 0 | 105 | 0 | 12.80 | 12.80 |
| 2 | 10 | 93 (80) | 1.07 (1.09) | 11.23 (9.72) | 12.3 (10.31) |
| 2 | 25 | 77 (84) | 2.71 (2.81) | 9.23 (10.39) | 11.94 (13.20) |
| 2 | 50 | 52 (77) | 5.54 (5.62) | 6.27 (9.49) | 11.81 (15.11) |
| 2 | 75 | 27 (67) | 8.61 (8.61) | 3.25 (8.25) | 11.86 (16.86) |
| 2 | 100 | 2 (61) | 11.37 (11.47) | 0.34 (7.60) | 11.71 (19.07) |
| 3 | 0 | 75 | 0 | 11.61 | 11.61 |
| 3 | 10 | 64 (57) | 1.41 (1.42) | 9.76 (8.66) | 11.17 (10.08) |
| 3 | 25 | 49 (57) | 3.61 (3.52) | 7.54 (8.63) | 11.51 (12.15) |
| 3 | 50 | 24 (51) | 7.20 (7.15) | 3.66 (7.76) | 10.86 (14.91) |
| 3 | 75 | 3 (44) | 11.00 (10.95) | 0.49 (6.70) | 11.49 (17.65) |
| 3 | 100 | 2 (36) | 14.75 (14.95) | 0.49 (5.62) | 15.24 (20.57) |
| 4 | 0 | 78 | 0 | 14.22 | 14.22 |
| 5 | 0 | 62 | 0 | 13.47 | 13.47 |
| 10 | 0 | 49 | 0 | 18.46 | 18.46 |
| 100 | 0 | 17 | 0 | 56.58 | 56.58 |